\documentclass{fic-l}
\input{diagrams}

\newtheorem{theorem}{Theorem}[section]

\newtheorem{proposition}[theorem]{Proposition}
\newtheorem{corollary}[theorem]{Corollary}

\theoremstyle{definition}
\newtheorem{definition}[theorem]{Definition}
\newtheorem{example}[theorem]{Example}

\theoremstyle{remark}
\newtheorem{remark}[theorem]{Remark}

\numberwithin{equation}{section}

\newcommand{\thlabel}[1]{\label{th:#1}}
\newcommand{\thref}[1]{Theorem~\ref{th:#1}}
\newcommand{\selabel}[1]{\label{se:#1}}
\newcommand{\seref}[1]{Section~\ref{se:#1}}

\newcommand{\prlabel}[1]{\label{pr:#1}}
\newcommand{\prref}[1]{Proposition~\ref{pr:#1}}
\newcommand{\colabel}[1]{\label{co:#1}}
\newcommand{\coref}[1]{Corollary~\ref{co:#1}}
\newcommand{\relabel}[1]{\label{re:#1}}

\newcommand{\exlabel}[1]{\label{ex:#1}}
\newcommand{\exref}[1]{Example~\ref{ex:#1}}
\newcommand{\delabel}[1]{\label{de:#1}}
\newcommand{\deref}[1]{Definition~\ref{de:#1}}
\newcommand{\eqlabel}[1]{\label{eq:#1}}
\newcommand{\equref}[1]{(\ref{eq:#1})}

\newcommand{\Hom}{{\rm Hom}}
\newcommand{\End}{{\rm End}}

\newcommand{\can}{{\rm can}}
\def\Ab{\underline{\underline{\rm Ab}}}
\def\lan{\langle}
\def\ran{\rangle}
\def\ot{\otimes}

\def\ZZ{{\mathbb Z}}
\def\QQ{{\mathbb Q}}

\newcommand{\Cc}{\mathcal{C}}
\newcommand{\Dd}{\mathcal{D}}

\newcommand{\Mm}{\mathcal{M}}

\def\text#1{{\rm {\rm #1}}}

\def\ol{\overline}
\def\ul{\underline}
\def\dul#1{\underline{\underline{#1}}}

\def\*C{{{}^*\mathcal C}}
\def\leftact{\hbox{$\rightharpoonup$}}

\begin{document}

\title{Galois corings from the descent theory point of view}

\author{S. Caenepeel}
\address{Faculty of Applied Sciences\\
Vrije Universiteit Brussel, VUB\\
Pleinlaan 2\\
B-1050 Brussels, Belgium}

%    General info
%\date{July 2, 1991}


\begin{abstract}
We introduce Galois corings, and give a survey of properties
that have been obtained so far. The Definition is motivated
using descent theory, and we show that classical Galois theory,
Hopf-Galois theory and coalgebra Galois theory can be obtained
as a special case.
\end{abstract}
\maketitle

\section*{Introduction}\selabel{0}
Galois descent theory \cite{Serre} has many applications in several branches
of mathematics, such as number theory, commutative algebra and
algebraic geometry; to name such one example, it is an essential
tool in computing the Brauer group of a field. In the literature,
several generalizations have appeared. Galois theory of commutative
rings has been studied by Auslander and Goldman \cite{AG}
and by Chase, Harrison and Rosenberg \cite{CHR}, see also \cite{DI}.
The group action can be replaced by a Hopf algebra (co)action,
leading to Hopf-Galois theory, see \cite{CS} (in the case where
the Hopf algebra is finitely generated projective), and \cite{DT},
\cite{Schneider90} in the general case. More recently,
coalgebra Galois extensions were introduced by 
Brzezi\'nski and Hajac \cite{BrzezinskiH}. It became clear recently
that a nice unification of all these theories can be formulated
using the language of corings. Let us briefly sketch the history.

During the nineties, several unifications of the various kinds
of Hopf modules that had appeared in the literature have been
proposed. Doi \cite{Doi92} and Koppinen \cite{Koppinen95}
introduced Doi-Hopf modules.
A more general concept, entwined modules, was proposed by
Brzezi\'nski and Majid \cite{BrzezinskiM}. B\"ohm 
introduced Doi-Hopf modules over a weak bialgebra (\cite{Bohm}),
and the author and De Groot proposed weak entwined modules
\cite{CaenepeelDG00}. Takeuchi \cite{Tak} observed that
all types of modules can be viewed as comodules over a coring,
a concept that was already introduced by Sweedler \cite{Sweedler65},
but then more or less forgotten, at least by Hopf algebra theorists;
the idea was further investigated by Brzezi\'nski \cite{Br}.
He generalized several properties that had been studied in
special cases to the situation where one works over a general
coring, such as separability and Frobenius type properties, and it
turned out that computations sometimes become amazingly simple if
one uses the language of corings, indicating that this is really
the right way to look at the problem.  Brzezi\'nski also
introduces the notion of Galois coring: to a ring extension
$i:\ B\to A$, one
can associate the so-called canonical coring; a morphism from
the canonical coring to another coring $\Cc$ is determined completely
by a grouplike element $x$; if this morphism is an isomorphism,
then we say that $(\Cc,x)$ is a Galois coring.

The canonical coring leads to an elegant formulation of descent
theory: the category of descent data associated to the
extension  $i:\ B\to A$ is nothing else then the category of
comodules over the canonical coring. This is no surprise:
to an $A$-coring, we can associate a comonad on the category of 
$A$-modules, and the canonical coring is exactly the comonad
associated to the adjoint pair of functors, given by induction and
restriction of scalars. Thus, if a coring is isomorphic to the
canonical coring, and if the induction functor is
comonadic, then the category of descent data is isomorphic to
the category of comodules over this coring, and equivalent to
the category of $B$-modules. This unifies all the versions of
descent theory that we mentioned at the beginning of this
introduction.

In this paper, we present a survey of properties of Galois corings
that have been obtained so far. We have organized it as follows:
in \seref{1}, we recall definition, basic properties and examples
of comodules over corings; in \seref{2}, we explain how descent
theory can be formulated using the canonical coring. We included
a full proof of \prref{2.3}, which is the noncommutative version
of the fact that the induction functor is comonadic if and only
if the ring morphism is pure as a map of modules. In \seref{3},
we introduce Galois corings, and discuss some properties, taken
from \cite{Br} and \cite{W3}. In \seref{4}, Morita theory is
applied to find some equivalent properties for a 
progenerator coring to be Galois; in fact, \thref{4.7} is a
new result, and generalizes results of Chase and Sweedler
\cite{CS}. In \seref{5}, we look at special cases,  and we show
how to recover the ``old" Galois theories. In \seref{6}, we
present a recent generalization, due to El Kaoutit and
G\'omez Torrecillas \cite{Kaoutit}.

\section{Corings}\selabel{1}
Let $A$ be a ring (with unit). The category of $A$-bimodules is
a braided monoidal category, and an $A$-coring is by definition a
coalgebra in the category of $A$-bimodules. Thus an $A$-coring
is a triple $\Cc=(\Cc,\Delta_\Cc,\varepsilon_\Cc)$, where
\begin{itemize}
\item $\Cc$ is an $A$-bimodule;
\item $\Delta_\Cc:\ \Cc\to \Cc\ot_A\Cc$ is an $A$-bimodule map;
\item $\varepsilon_\Cc:\ \Cc\to A$ is an $A$-bimodule map
\end{itemize}
such that
\begin{equation}\eqlabel{1.1}
(\Delta_\Cc\ot_A I_\Cc)\circ \Delta_\Cc=
(I_\Cc\ot_A\Delta_\Cc )\circ \Delta_\Cc,
\end{equation}
and
\begin{equation}\eqlabel{1.2}
(I_{\Cc}\ot_A\varepsilon_{\Cc})\circ \Delta_\Cc=
(\varepsilon_{\Cc}\ot_AI_{\Cc})\circ \Delta_\Cc=I_\Cc.
\end{equation}
Sometimes corings are considered as coalgebras over noncommutative
rings. This point of view is not entirely correct: a coalgebra over
a commutative ring $k$ is a $k$-coring, but not conversely: it could
be that the left and and right action of $k$ on the coring are different.

The Sweedler-Heyneman notation is also used for a coring $\Cc$,
namely
$$\Delta_\Cc(c)=c_{(1)}\ot_A c_{(2)},$$
where the summation is implicitely understood.  \equref{1.2} can
then be written as
$$\varepsilon_\Cc(c_{(1)})c_{(2)}=c_{(1)}\varepsilon_{\Cc}(c_{(2)})=c.$$
This formula looks like the corresponding formula for usual coalgebras.
Notice however that the order matters in the above formula, since
$\varepsilon_{\Cc}$ now takes values in $A$ which is noncommutative
in general. Even worse, the expression 
$c_{(2)}\varepsilon_{\Cc}(c_{(1)})$ makes no sense at all, since we
have no well-defined switch map $\Cc\ot_A\Cc\to \Cc\ot_A\Cc$.\\
A morphism between two corings $\Cc$ and $\Dd$ is an $A$-bimodule
map $f:\ \Cc\to \Dd$ such that
$$\Delta_\Dd(f(c))=f(c_{(1)})\ot_A  f(c_{(2)})~~{\rm and}~~
\varepsilon_\Dd(f(c))=\varepsilon_\Cc(c),$$
for all $c\in \Cc$.
A right $\Cc$-comodule $M=(M,\rho )$ consists of a right $A$-module
$M$ together with a right $A$-linear map $\rho :\ M\to M\ot_A\Cc$
such that:
\begin{equation}\eqlabel{1.3}
(\rho \ot_A I_\Cc)\circ \rho =(I_M\ot_A\Delta_\Cc)\circ \rho, 
\end{equation}
and
\begin{equation}\eqlabel{1.4}
(I_M\ot_A\varepsilon_\Cc)\circ \rho =I_M.
\end{equation}
We then say that $\Cc$ coacts from the right on $M$.
Left $\Cc$-comodules and $\Cc$-bicomodules can be defined in a similar
way. We use the Sweedler-Heyneman notation also for comodules:
$$\rho (m)=m_{[0]}\ot_A m_{[1]}.$$
\equref{1.4} then takes the form
$m_{[0]}\varepsilon_\Cc( m_{[1]})=m$. A right $A$-linear map
$f:\ M\to N$ between two right $\Cc$-comodules $M$ and $N$ is called
right $\Cc$-colinear if $\rho (f(m))=f(m_{[0]})\ot m_{[1]}$,
for all $m\in M$.

Corings were already considered by Sweedler in \cite{Sweedler65}.
The interest in corings was revived after a mathematical review
written by Takeuchi \cite{Tak}, in which he observed that
entwined modules can be considered as comodules over a coring.
This will be discussed in the examples below.

\begin{example}\exlabel{1.2}
As we already mentioned, if $A$ is a commutative ring, then
an $A$-coalgebra is also an $A$-coring.
\end{example}

\begin{example}\exlabel{1.3}
Let $i:\ B\to A$ be a ring morphism; then $\Dd=A\ot_B A$
is an $A$-coring. We define
$$\Delta_\Dd:\ \Dd\to \Dd\ot_A\Dd\cong A\ot_B A\ot_B A$$
and
$$\varepsilon_\Dd:\ \Dd=A\ot_B A\to A$$
by
$$\Delta_\Dd(a\ot_B b)=(a\ot_B 1_A)\ot_A(1_A\ot_B b)\cong
a\ot_B 1_A\ot_B b$$
and
$$\varepsilon_{\Dd}(a\ot_B b)=ab.$$
Then $\Dd=(\Dd,\Delta_\Dd,\varepsilon_{\Dd})$ is an $A$-coring. It 
is called the canonical coring associated to the ring morphism $i$.
We will see in the next section that this coring is crucial in 
descent theory.
\end{example}

\begin{example}\exlabel{1.4}
Let $k$ be a commutative ring, $G$ a finite group, and
$A$ a $G$-module algebra. Let
$\Cc=\oplus_{\sigma\in G} Av_{\sigma}$
be the left free $A$-module with basis indexed by $G$,
and let $p_\sigma:\  \Cc\to A$ be the projection onto the free component
$Av_\sigma$.
We make
$\Cc$ into a right $A$-module by putting
$$v_{\sigma}a=\sigma(a)v_{\sigma}.$$
A comultiplication and counit on $\Cc$ are defined by putting
$$\Delta_\Cc(av_\sigma)=\sum_{\tau\in G} av_{\tau}\ot_A v_{\tau^{-1}\sigma}
~~{\rm and}~~\varepsilon_\Cc=p_e,$$
where $e$ is the unit element of $G$. It is straightforward to
verify that $\Cc$ is an $A$-coring. Notice that, in the case
where $A$ is commutative, we have an example of an $A$-coring,
which is not an $A$-coalgebra, since the left and right $A$-action
on $\Cc$ do not coincide.

Let us give a description of the right $\Cc$-comodules. Assume
that $M=(M,\rho )$ is a right $\Cc$-comodule. For every $m\in M$
and $\sigma\in G$, let $\ol{\sigma}(m)=m_{\sigma}=I_M\ot_A p_{\sigma})(\rho (m))$.
Then we have
$$\rho (m)=\sum_{\sigma\in G} m_\sigma\ot_A v_\sigma.$$
$\ol{e}$ is the identity, since
$m=(I_M\ot_A\varepsilon_\Cc)\circ\rho (m)=m_e$.
Using the coassociativity of the comultiplication, we find
$$
\sum_{\sigma\in G}\rho (m_\sigma)\ot v_\sigma=
\sum_{\sigma,\tau\in G}m_\sigma\ot_Av_\tau\ot_A v_{\tau^{-1}\sigma}
=\sum_{\rho,\tau\in G}m_{\tau\rho}\ot_Av_\tau\ot_A v_{\rho},
$$
hence $\rho (m_\sigma)=\sum_{\tau\in G} m_{\tau\sigma}\ot_A v_\tau$,
and $\ol{\tau}(\ol{\sigma}(m))=m_{\tau\sigma}=\ol{\tau\sigma}(m)$,
so $G$ acts as a group of $k$-automorphisms on $M$. Moreover, since
$\rho $ is right $A$-linear, we have that
$$
\rho (ma)=\sum_{\sigma\in G}\ol{\sigma}(ma)\ot_Av_\sigma
=\sum_{\sigma\in G}\ol{\sigma}(m)\ot_Av_\sigma a=
\sum_{\sigma\in G}\ol{\sigma}(m)\sigma(a)\ot_Av_\sigma
$$
so $\ol{\sigma}$ is $A$-semilinear (cf. \cite[p. 55]{KO}):
$\ol{\sigma}(ma)=\ol{\sigma}(m)\sigma(a)$, for all $m\in M$ and $a\in A$.
Conversely, if $G$ acts as a group of right $A$-semilinear
automorphims on $M$, then the formula
$$\rho (m)=\sum_{\sigma\in G} \ol{\sigma}(m)\ot_Av_\sigma$$
defines a right $\Cc$-comodule structure on $\Mm$.
\end{example}

\begin{example}\exlabel{1.5}
Now let $k$ be a commutative ring, $G$ an arbitrary group,
and $A$ a $G$-graded $k$-algebra. Again let $\Cc$ be the free
left $A$-module with basis indexed by $G$:
$$\Cc=\oplus_{\sigma\in G} Au_{\sigma}$$
Right $A$-action, comultiplication and counit are now defined by
$$u_\sigma a=\sum_{\tau\in G} a_\tau u_{\sigma\tau}~~;~~
\Delta_\Cc(u_\sigma)=u_{\sigma}\ot_A u_{\sigma}~~;~~
\varepsilon_\Cc(u_\sigma)=1.$$
$\Cc$ is an $A$-coring; let $M=(M,\rho )$ be a right $A$-comodule,
and let $M_\sigma=\{m\in M~|~\rho(m)=m\ot_A u_\sigma\}$.
It is then clear that $M_\sigma\cap M_\tau=\{0\}$ if $\sigma\neq \tau$.
For any $m\in M$, we can write in a unique way:
$$\rho (m)=\sum_{\sigma\in G} m_{\sigma}\ot_A u_{\sigma}.$$
Using the coassociativity, we find that $m_{\sigma}\in M_\sigma$,
and using the counit property, we find that $m=\sum_{\sigma} m_\sigma$.
So $M=\oplus_{\sigma\in G} M_\sigma$. Finally, if $m\in M_\sigma$
and $a\in A_\tau$, then it follows from the right $A$-linearity
of $\rho $ that
$$\rho (ma)=(m\ot_A u_\sigma)a=ma\ot_A u_{\sigma\tau},$$
so $ma\in M_{\sigma\tau}$, and $M_\sigma A_\tau\subset M_{\sigma\tau}$,
and $M$ is a right $G$-graded $A$-module. Conversely, every
right $G$-graded $A$-module can be made into a right $\Cc$-comodule.
\end{example}

\begin{example}\exlabel{1.6}
Let $H$ be a bialgebra over a commutative ring $k$, 
and $A$ a right $H$-comodule algebra. Now take
$\Cc=A\ot H$, with $A$-bimodule structure
$$a'(b\ot h)a=a'ba_{[0]}\ot ha_{[1]}.$$
Now identify $(A\ot H)\ot_A(A\ot H)\cong A\ot H\ot H$,
and define the comultiplication and counit on $\Cc$, by putting
$\Delta_\Cc=I_A\ot \Delta_H$ and $\varepsilon_\Cc=I_A\ot \varepsilon_H$.
Then $\Cc$ is an $A$-coring.
The category $\Mm^\Cc$ is isomorphic to the category of relative Hopf
modules. These are $k$-modules $M$ with a right $A$-action and a right
$H$-coaction $\rho$, such that
$$\rho(ma)=m_{[0]}a_{[0]}\ot_Am_{[1]}a_{[1]}$$
for all $m\in M$ and $a\in A$.
\end{example}

\begin{example}\exlabel{1.7}
Let $k$ be a commutative ring, $A$ a $k$-algebra, and $C$ a $k$-coalgebra,
and consider a $k$-linear map
$\psi:\ C\ot A\to A\ot C$.
We use the following Sweedler type notation, where the summation is
implicitely understood:
$$\psi(c\ot a)=a_{\psi}\ot c^{\psi}=a_{\Psi}\ot c^{\Psi}.$$
$(A,C,\psi)$ is called a (right-right) entwining structure if
the four following conditions are satisfied:
\begin{eqnarray}
&&  (ab)_{\psi}\ot c^{\psi}=  a_{\psi}b_{\Psi}\ot
c^{\psi\Psi};\eqlabel{1.7.1}\\
&&(1_A)_{\psi}\ot c^{\psi}=1_A\ot c;\eqlabel{1.7.2}\\
&&  a_{\psi}\ot \Delta_C(c^{\psi})=
  a_{\psi\Psi}\ot c_{(1)}^{\Psi}\ot c_{(2)}^{\psi};
\eqlabel{1.7.3}\\
&&\varepsilon_C(c^{\psi})  a_{\psi} =
\varepsilon_C(c)a.\eqlabel{1.7.4}
\end{eqnarray}
Let $\Cc=A\ot C$ as a $k$-module, with $A$-bimodule structure
$$a'(b\ot c)a=a'ba_{\psi}\ot c^{\psi}.$$
Comultiplication and counit on $A\ot C$ are defined as in
\exref{1.6}. $\Cc$ is a coring, and the category $\Mm^\Cc$ is
isomorphic to the category $\Mm(\psi)_A^C$
of entwined modules. These are $k$-modules
$M$ with a right $A$-action and a right $C$-coaction $\rho$ such that
$$\rho(ma)=m_{[0]}a_{\psi}\ot_Am_{[1]}^\psi,$$
for all $m\in M$ and $a\in A$.
\end{example}

Actually Examples \ref{ex:1.4}, \ref{ex:1.5} and \ref{ex:1.6}
are special cases of \exref{1.7}
\begin{itemize}
\item \exref{1.4}: take
$C=(kG)^*=\oplus_{g\in G} kv_\sigma$, the dual of the group ring
$kG$, and
$\psi(v_\sigma\ot a)=\sigma(a)\ot v_\sigma$;
\item \exref{1.5}: take $C=kG=\oplus_{g\in G}ku_\sigma$, the group ring,
and $\psi(u_\sigma\ot a)=\sum_{\tau\in G} a_\tau\ot u_{\sigma\tau}$;
\item \exref{1.6}: take $C=H$, and $\psi(h\ot a)=a_{[0]}\ot ha_{[1]}$.
\end{itemize}

If $\Cc$ is an $A$-coring, then its left dual
${\*C}={}_A\Hom(\Cc,A)$
is a ring, with (associative) multiplication
given by the formula
\begin{equation}\eqlabel{1.8.1}
f\# g=g\circ (I_\Cc\ot_A f)\circ \Delta_\Cc
~~{\rm or}~~
(f\# g)(c)=g(c_{(1)}f(c_{(2)})),
\end{equation}
for all left $A$-linear $f,g:\ \Cc\to A$ and $c\in \Cc$.
The unit is $\varepsilon_\Cc$. We have 
a ring homomorphism
$i:\ A\to {\*C},~~i(a)(c)=\varepsilon_\Cc(c)a$.
We easily compute that
\begin{equation}\eqlabel{1.8.2}
(i(a)\#f)(c)=f(ca)~~{\rm and}~~(f\# i(a))(c)=f(c)a,
\end{equation}
for all $f\in {\*C}$, $a\in A$ and $c\in \Cc$.
We have a functor
$F:\ \Mm^\Cc\to \Mm_{{\*C}}$, where
$F(M)=M$ as a right $A$-module, with right ${\*C}$-action
given by
$m\cdot f=m_{[0]}f(m_{[1]})$, for all $m\in M$, $f\in \*C$.
If $\Cc$ is finitely generated and projective as a left
$A$-module, then $F$ is an isomorphism of categories: given
a right ${\*C}$-action on $M$, we recover the right
$\Cc$-coaction by putting
$\rho(m)=\sum_j (m\cdot f_j)\ot_A c_j$,
where $\{(c_j,f_j)~|~j=1,\cdots, n\}$ is a finite dual basis of $\Cc$
as a left $A$-module.
${\*C}$ is a right $A$-module, by \equref{1.8.2}:
$(f\cdot a)(c)=f(c)a$,
and we can consider the double dual $({\*C})^*=\Hom_A({\*C},A)$.
We have a canonical morphism
$i:\ \Cc\to ({\*C})^*,~~i(c)(f)=f(c)$, and
we call $\Cc$ reflexive (as a left $A$-module) if $i$ is
an isomorphism. If $\Cc$ is finitely generated projective as
a left $A$-module, then $\Cc$ is reflexive. For any $\varphi\in
({\*C})^*$, we then have that
$\varphi=i(\sum_j \varphi(f_j)c_j)$.

\subsubsection*{Corings with a grouplike element}
Let $\Cc$ be an $A$-coring, and suppose that $\Cc$ coacts on $A$.
Then we have a map $\rho :\ A\to A\ot_A\Cc\cong \Cc$. The fact that
$\rho $ is right $A$-linear implies that $\rho $ is completely
determined by $\rho (1_A)=x$: $\rho (a)=xa$. 
The coassociativity of the coaction
yields that $\Delta_{\Cc}(x)=x\ot_A x$ and the counit property gives
us that $\varepsilon_{\Cc}(x)=1_A$. We say that $x$ is a {\it grouplike
element} of $\Cc$ and we denote $G(\Cc)$ for the set of all grouplike
elements of $\Cc$. If $x\in G(\Cc)$ is grouplike, then the associated
$\Cc$-coaction on $A$ is given by $\rho (a)=xa$.

If $x\in G(\Cc)$, then we call $(\Cc,x)$ a {\it coring with a fixed
grouplike element}. For $M\in \Mm^{\Cc}$, we call
$$M^{{\rm co}\Cc}=\{m\in M~|~\rho (m)=m\ot_A x\}$$
the submodule of coinvariants of $M$; note that this definition depends
on the choice of the grouplike element. Also observe that
$$A^{{\rm co}\Cc}=\{b\in A~|~bx=xb\}$$
is a subring of $A$. Let $i:\ B\to A$ be a ring morphism.
$i$ factorizes through $A^{{\rm co}\Cc}$ if and only if
$$x\in G(\Cc)^B=\{x\in G(\Cc)~|~xb=bx,~{\rm for~all~}b\in B\}.$$ 
We then have
two pairs of adjoint functors
$(F,G)$ and $(F',G')$, respectively between the categories $\Mm_B$ and
$\Mm^\Cc$ and the categories ${}_B\Mm$ and
${}^\Cc\Mm$. For 
$N\in \Mm_B$ and $M\in \Mm^\Cc$,
$$F(N)=N\ot_B A~~{\rm and}~~G(M)=M^{{\rm co}\Cc}.$$
The unit and counit of the adjunction are
$$\nu_N:\ N\to (N\ot_BA)^{{\rm co}\Cc},~~\nu_N(n)=n\ot_B1;$$
$$\zeta_M:\ M^{{\rm co}\Cc}\ot_B A\to M,~~
\zeta_M(m\ot_B a)=ma.$$
The other adjunction is defined in a similar way. We want to
discuss when $(F,G)$ and $(F',G')$ are category equivalences.
In \seref{2}, we will do this for the canonical coring
associated to a ring morphism; we will study the general case
in \seref{3}.

\section{The canonical coring and descent theory}\selabel{2}
Let $i:\ B\to A$ be a ring morphism. The problem of descent theory
is the following: suppose that we have a right $A$-module $M$.
When do we have a right $B$-module $N$ such that $N=M\ot_BA$?
The same problem can be stated for modules with an additional structure,
for example algebras. In the case where $A$ and $B$ are commutative,
this problem has been discussed in a purely algebraic context
in \cite{KO}. In fact the results in \cite{KO} are the affine versions
of Grothendieck's descent theory for schemes, see \cite{Gr}.
In the situation where $A$ and $B$ are arbitrary, descent theory has
been discussed by Cipolla \cite{Cipolla}, and, more recently,
by Nuss \cite{Nuss}. For a purely categorical treatment of the problem,
making use of monads, we refer to \cite{Borceux}. Here we will
show that the results in \cite{KO} and \cite{Cipolla} can be
restated elegantly in terms of comodules over the canonical coring.

Let $\Dd=A\ot_B A$ be the canonical coring associated to the ring morphism
$i:\ B\to A$, and let $M=(M,\rho )$ be a right $\Dd$-comodule.
We will identify $M\ot_A\Dd\cong M\ot_B A$ using the natural isomorphism.
The coassociativity and the counit property then take the form
$$\rho (m_{[0]})\ot m_{[1]}=m_{[0]}\ot_B 1_A\ot_B m_{[1]}~~{\rm and}~~
m_{[0]}m_{[1]}=m.$$
$1_A\ot_B1_A$ is a grouplike element of $\Dd$. As we have seen
at the end of \seref{1}, we have two pairs of adjoint functors,
respectively between $\Mm_B$ and $\Mm^\Dd$, and
${}_B\Mm$ and ${}^\Dd\Mm$, which we will denote by $(K,R)$ and $(K',R')$.
The unit and counit of the adjunction
will be denoted by $\eta$ and $\varepsilon$.
$K$ is called
the comparison functor.
If $(K,R)$ is an equivalence of categories, then the ``descent problem"
is solved: $M\in \Mm_A$ is isomorphic to some $N\ot_B A$ if and only
if we can define a right $\Dd$-coaction on $M$.

Recall that a morphism of left $B$-modules $f:\ M\to M'$
is called pure if and only if $f_N=I_N\ot_B f:\ N\ot_B M\to N\ot_B M'$
is monic, for every $N\in \Mm_B$. $i:\ B\to A$ is pure as a morphism of left $B$-modules if
and only if
 $\eta_N$ being injective, for all $N\in \Mm_B$, since
 $\eta_N$ factorizes through $i_N$.

\begin{proposition}\prlabel{2.1}
The comparison functor $K$ is fully faithful if and only if
$i:\ B\to A$ is pure as a morphism of left $B$-modules.
\end{proposition}

\begin{proof}
The comparison functor $K$ is fully faithful if and only if
$\eta_N$ is bijective, for all $N\in\Mm_B$.  From the above observation,
it follows that it suffices to show that left purity of $i$ implies that
$\eta_N$ is surjective. Since $\eta_N$ is injective, we have
that $N\subset (N\ot_BA)^{{\rm co}\Dd}\subset N\ot_B A$.
Take $q=\sum_i n_i\ot_B a_i\in (N\ot_BA)^{{\rm co}\Dd}$. Then
\begin{equation}\eqlabel{2.1.1}
\rho(\sum_i n_i\ot_B a_i)=\sum_i n_i\ot_B a_i\ot_B 1=\sum_i n_i\ot_B1\ot_B a_i.
\end{equation}
Consider the right $B$-module $P=(P\ot_B A)/N$, and let
$\pi:\ N\ot_B A\to P$ be the canonical projection. Applying $\pi\ot_BI_A$ to
\equref{2.1.1}, we obtain
$$\pi(q)\ot_B 1=\sum_i \pi(n_i\ot_B 1)\ot_Ba_i=0~{\rm in}~P\ot_B A,$$
hence $\pi(q)=0$, since $i_P$ is an injection. This means that
$q\in N$, which is exactly what we needed.
\end{proof}

We also have an easy characterization of the fact that $R$ is fully
faithful.

\begin{proposition}\prlabel{2.2}
The right adjoint $R$ of the comparison functor $K$ is fully faithful if and only
if
$\bullet\ot_B A$ preserves the equalizer of $\rho $ and $i_M$,
for every $(M,\rho )\in \Mm^\Dd$. In particular, if $A$ is flat
as a left $B$-module, then $R$ is fully faithful.
\end{proposition}

\begin{proof}
$M^{{\rm co}\Dd}$ is the equalizer of the maps
$$0\rTo M^{{\rm co}\Dd}\rTo^{j} M 
\pile{\rTo^{\rho}\\ \rTo_{i_M}} M\ot_B A.$$
First assume that $M^{{\rm co}\Dd}\ot_B A$ is the equalizer
\begin{equation}\eqlabel{2.2.1}
0\rTo M^{{\rm co}\Dd}\ot_B A\rTo^{j\ot_B I_A} M \ot_B A
\pile{\rTo^{\rho\ot_BI_A}\\ \rTo_{i_M\ot_BI_A}} M\ot_B A\ot_B A.
\end{equation}
From the coassociativity of $\rho$, it now follows that
$\rho (m)\in M^{{\rm co}\Dd}\ot_B A\subset M\ot_B A\cong
M\ot_A(A\ot_B A)$, for all $m\in M$, and we have a map
$\rho:\ M\to M^{{\rm co}\Dd}\ot_B A$. From the counit property,
it follows that $\varepsilon_M\circ \rho=I_M$. For $m\in 
M^{{\rm co}\Dd}$ and $a\in A$, we have
$$\rho(\varepsilon_M(m\ot_B a))=\rho(ma)=\rho(m)a=
(m\ot_B 1)a=m\ot_B a.$$
Thus the counit $\varepsilon_M$ has an inverse, for all $M$, and
$R$ is fully faithful.

Conversely, assume that $\varepsilon_M$ is bijective. Take
$\sum_i m_i\ot_Ba_i\in M^{{\rm co}\Dd}\ot_B A$, and put
$m=\sum_i m_ia_i\in M$. Then $\rho (m)=m_{[0]}\ot_B m_{[1]}=
\sum_i m_i\ot_Ba_i\in M\ot_B A$. Consequently, if 
$\sum_i m_i\ot_Ba_i=0\in M\ot_B A$, then $m=m_{[0]}m_{[1]}=0$,
so $\sum_i m_i\ot_Ba_i=0\in M^{{\rm co}\Dd}\ot_B A$, and we have
shown that the canonical map $M^{{\rm co}\Dd}\ot_B A\to
M\ot_B A$ is injective.
\end{proof}

If $A$ and $B$ are commutative, and
$i:\ B\to A$ is pure as a morphism of $B$-modules, then
$\bullet\ot_B A$ preserves the equalizer of $\rho $ and $i_M$,
for every $(M,\rho )\in \Mm^\Dd$, and therefore $R$ is fully faithful.
This result is due to Joyal and Tierney (unpublished); an
elementary proof was given recently by Mesablishvili \cite{Mesablishvili}.
We will now adapt Mesablishvili's proof to the noncommutative
situation. In view of \prref{2.1}, one would expect that a
sufficient condition for the fully faithfulness of $R$ is the
fact that $i$ is pure as a morphism of left $B$-modules.
It came as a surprise to the author that we need right purity
instead of left purity.

We consider the contravariant functor
$C=\Hom_\ZZ(\bullet,\QQ/\ZZ):\ \Ab\to \Ab.$
$\QQ/\ZZ$ is an injective cogenerator of $\Ab$, and therefore
$C$ is exact and reflects isomorphisms. If $B$ is a ring, then
$C$ induces functors
$$C:\ \Mm_B\to {}_B\Mm~~{\rm and}~~{}_B\Mm\to \Mm_B.$$
For example, if $M\in \Mm_B$, then $C(M)$ is a left $B$-module,
by putting $(b\cdot f)(m)=f(mb)$.
For $M\in \Mm_B$ and $P\in{}_BM$, we have the following isomorphisms,
natural in $M$ and $P$:
\begin{equation}\eqlabel{2.3.1}
\Hom_B(M,C(P))\cong {}_B\Hom(P,C(M))\cong C(M\ot_B P)
\end{equation}
If $P\in {}_B\Mm_B$, then $C(P)\in {}_B\Mm_B$, and the above isomorphisms
are isomorphisms of left $B$-modules.

\begin{proposition}\prlabel{2.3}
Let $i:\ B\to A$ be a ring morphism, and assume that $i$ is pure
as a morphism of right $B$-modules. Then the adjoint $R$ of
the comparison functor is fully faithful.
\end{proposition}

\begin{proof}
We have to show that \equref{2.2.1} is exact, for all
$(M,\rho)\in \Mm^\Dd$.
If $i$ is pure in $\Mm_B$, then 
$i_{C(B)}:\ B\ot_BC(B)\to A\ot_B C(B)$ is a
monomorphism in $\Mm_B$, hence
$$C(i_{C(B)}):\ C(A\ot_B C(B))\to C(B\ot_B C(B))$$
is an epimorphism in ${}_B\Mm$. Applying \equref{2.3.1}, we find
that
$$C(i)\circ\bullet :\ {}_B\Hom(C(B),C(A))\to {}_B\Hom(C(B),C(B))$$
is also an epimorphism. This implies that $C(i):\
C(A)\to C(B)$ is a split epimorphism in ${}_B\Mm$, and then it follows
that for every $M\in \Mm_B$,
$$C(i)\circ\bullet:\ \Hom_B(M,C(A))\to \Hom_B(M,C(B))$$
is a split epimorphism in ${}_B\Mm$. Applying \equref{2.3.1}
again, we find that
$$C(i_M):\ C(M\ot_BA)\to C(M)$$
is a split epimorphism in ${}_B\Mm$.\\
In $\Mm_B$, we have the following commutative diagram with
exact rows:
$$\begin{diagram}
0&\rTo^{}& M^{{\rm co}\Dd}&\rTo^j&M&
\pile{\rTo^{\rho }\\ \rTo_{i_M}}&M\ot_BA\\
&&\dTo_{j}&&\dTo_{i_M}&&\dTo_{i_{M\ot_BA}}\\
0&\rTo^{}&M&\rTo^{\rho }&M\ot_BA&
\pile{\rTo^{\rho \ot I_A}\\ \rTo_{i_M\ot I_A}}&M\ot_BA\ot_BA
\end{diagram}$$
Applying the functor $C$ to this diagram, we obtain the following
commutative diagram with exact rows. 
We also know that
$C(i_M)$ and $C(i_{M\ot_BA})$ have right inverses $h$ and $h'$.
$$\begin{diagram}
C(M\ot_BA\ot_BA)&
\pile{\rTo^{C(\rho \ot I_A)}\\ \rTo_{C(i_M\ot I_A)}}&
C(M\ot_BA)&\rTo^{C(\rho )}&C(M)&\rTo^{}&0\\
\uTo^{h'}\dTo_{C(i_{M\ot_BA})}&&\uTo^{h}\dTo_{C(i_M)}
&&\uTo^{k}\dTo_{C(j)}\\
C(M\ot_BA)&
\pile{\rTo^{C(\rho )}\\ \rTo_{C(i_M)}}&
C(M)&\rTo^{C(j)}&C(M^{{\rm co}\Dd})&\rTo^{}& 0
\end{diagram}$$
Diagram chasing leads to the existence of a right inverse $k$
of $C(j)$, such that $k\circ C(j)=C(\rho )\circ h$. But this
means that the bottom row of the above diagram is a split fork
in ${}_B\Mm$,
split by the morphisms
$$C(M\ot_BA)\lTo^{h}C(M)\lTo^{k}C(M^{{\rm co}\Dd})$$
(see \cite[p.149]{McLane} for the definition of a split fork).
Split forks are preserved by arbitrary functors, so applying
${}_B\Hom(A,\bullet)$, we obtain a split fork
in ${}_B\Mm$; using \equref{2.3.1}, we find that this split fork
is isomorphic to
$$C(M\ot_BA\ot_BA)
\pile{\rTo^{C(\rho \ot I_A)}\\ \rTo_{C(i_M\ot I_A)}}
C(M\ot_BA)\rTo^{C(j\ot I_A)}C(M^{{\rm co}\Dd}\ot_BA)$$
The functor $C$ is exact and reflects isomorphisms, hence it
also reflects coequalizers. It then follows that \equref{2.2.1}
is an equalizer in $\Mm_B$, as needed.
\end{proof}

The converse of \prref{2.3} is not true in general: the natural
inclusion $i:\  \ZZ\to \QQ$ is not pure in $\Mm_\ZZ$, but the functor
$R$ is fully faithful. Indeed, if $(M,\rho)\in \Mm^{\Dd}$, then
$M$ is a $\QQ$-vector space, and 
$\rho:\ M\to M\ot_{\ZZ} \QQ\cong M$ is the identity, $M^{{\rm co}\Dd}
=M$, and $\varepsilon_M:\ M^{{\rm co}\Dd}\ot_{\ZZ} \QQ\cong M\to M$
is also the identity.

It would be interesting to know if there exists a ring morphism $i:\ B\to A$
which is pure in ${}_B\Mm$, but not in $\Mm_B$, and such that
$(K,R)$ is an equivalence of categories.

Consider $K'=A\ot_B\bullet:\ {}_B\Mm\to {}^{\Dd}\Mm$ and
$R'={}^{{\rm co}\Dd}(\bullet):\ {}^{\Dd}\Mm\to {}_B\Mm$.
The next result is an immediate consequence of
Propositions \ref{pr:2.2} and \ref{pr:2.3} and their left
handed versions, and can be viewed as the noncommutative
version of the Joyal-Tierney Theorem.

\begin{theorem}\thlabel{2.4}
Let $i:\ B\to A$ be a morphism of rings. Then the following assertions
are equivalent.
\begin{enumerate}
\item $(K,R)$ and $(K',R')$ are equivalences of categories;
\item $K$ and $K'$ are fully faithful;
\item $i$ is pure in $\Mm_B$ and ${}_B\Mm$.
\end{enumerate}
\end{theorem}

We have seen in \prref{2.2} that $R$ is fully faithful if
$A$ is flat as a left $B$-module.

\begin{proposition}\prlabel{2.5}
Let $i:\ B\to A$ be a morphism of rings, and assume that $A$ is 
flat
as a left $B$-module. Then $(K,R)$ is an equivalence of categories
if and only if $A$ is faithfully flat as a left $B$-module.
\end{proposition}

\begin{proof}
First assume that $A$ is faithfully flat as a left $B$-module.
It follows from \prref{2.1} that it suffices to show that
$A$ is pure as a left $B$-module. For $N\in \Mm_B$, the map
$$f=i_N\ot_B I_A:\ N\ot_BA\to N\ot_BA\ot_BA$$
is injective: if $f(\sum_i n_i\ot_B a_i)=\sum_i n_i\ot_B1\ot_B a_i=0$,
then, multiplying the second and third tensor fact, we find
that $\sum_i n_i\ot_B a_i=0$. Since $A$ is faithfully flat as
a left $B$-module, it follows that $i_N$ is injective.\\
Conversely assume that $(K,R)$ is an equivalence of categories.
Then the functor $R$ is exact. Let
\begin{equation}\eqlabel{2.5.1}
0\to N'\to N\to N''\to 0
\end{equation}
be a sequence of right $B$-modules such that the sequence
$$0\to N'\ot_BA\to N\ot_BA\to N''\ot_BA\to 0$$
is exact. Applying the functor $R$ to the sequence, and using the
fact that $\eta$ is an isomorphism, we find that \equref{2.5.1} is
exact, so it follows that $A$ is faithfully flat as a left
$B$-module.
\end{proof}

The descent data that are considered in \cite{Cipolla} are nothing
else then comodules over the canonical coring (although the author
of \cite{Cipolla} was not aware of this). The descent data in \cite{KO}
are different, so let us indicate how to go from descent data to
comodules over the canonical coring.

Let $i:\ B\to A$ be a morphism of commutative rings. A {\it descent
datum} consists of a pair $(M,g)$,
with $M\in \Mm_A$, and $g:\ A\ot_B M\to M\ot_B A$ an 
$A\ot_BA$-module homorphism such that
\begin{equation}\eqlabel{2.6.1}
g_2=g_3\circ g_1:\ A\ot_BA\ot_B M\to A\ot_BM\ot_B A
\end{equation}
and
\begin{equation}\eqlabel{2.6.2}
\mu_M(g(1\ot_B m))=m,
\end{equation}
for all $m\in M$. Here $g_i$ is obtained by applying $I_A$ to the $i$-th
tensor position, and $g$ to the two other ones. It can be shown
that \equref{2.6.2} can be replaced by the condition that $g$ is
a bijection. A morphism of two
descent data $(M,g)$ and $(M',g')$ consists of an $A$-module homomorphism
$f:\ M\to M'$ such that 
$$(f\ot_B I_A)\circ g=g'\circ (I_A\ot_B f).$$
$\dul{\rm Desc}(A/B)$ will be the category of descent data.

\begin{proposition}\prlabel{2.6}
Let $i:\ B\to A$ be a morphism of commutative rings. We have
an isomorphism of categories
$$\dul{\rm Desc}(A/B)\cong {\Mm}^{A\ot_BA}$$
\end{proposition}

\begin{proof} (sketch)
For a right $\Dd$-comodule $(M,\rho )$, we define
$g:\ A\ot_B M\to M\ot_B A$ by $g(a\ot_B m)=m_{[0]}a\ot_B m_{[1]}$.
Then $(M,g)$ is a descent datum. Conversely, given a descent datum
$(M,g)$, the map $\rho :\ M\to M\ot_B A$, $\rho (m)=
g(1\ot_B m)$ makes $M$ into a right $\Dd$-comodule.
\end{proof}

\section{Galois corings}\selabel{3}
Let $A$ be a ring, $(\Cc,x)$ a coring with a fixed grouplike
element, and $i:\ B\to A^{{\rm co}\Cc}$ a ring morphism.
We have seen at the end of \seref{1} that we have two pairs
of adjoint functors $(F,G)$ and $(F',G')$. We also have
a morphism of corings
$${\rm can}:\ \Dd=A\ot_BA\to \Cc,~~{\rm can}(a\ot_B a')=axa'.$$

\begin{proposition}\prlabel{3.0}
With notation as above, we have the following results.
\begin{enumerate}
\item If $F$ is fully faithful, then $i:\ B\to A^{{\rm co}\Cc}$
is an isomorphism;
\item if $G$ is fully faithful, then ${\rm can}:\ \Dd=A\ot_BA\to \Cc$
is an isomorphism.
\end{enumerate}
\end{proposition}

\begin{proof}
1. $F$ is fully faithful if and only if $\nu$ is an isomorphism;
it then suffices to observe that $i=\nu_B$.\\
2. $G$ is fully faithful if and only if $\zeta$ is an isomorphism.
$\Cc\in\Mm^\Cc$, the right coaction being induced by the
comultiplication. The map $f:\ A\to \Cc^{{\rm co}\Cc}$,
$f(a)=ax$,
is an isomorphism of $(A,B)$-bimodules; the inverse of $f$ is
the restriction of $\varepsilon_\Cc$ to $\Cc^{{\rm co}\Cc}$.
Indeed, if $c\in \Cc^{{\rm co}\Cc}$, then $\Delta_\Cc(c)=
c\ot_A x$, hence $c=\varepsilon(c)x=f(\varepsilon(c))$. It 
follows that ${\rm can}=\zeta_\Cc\circ (f\ot_BI_A)$ is an
isomorphism.
\end{proof}

\prref{3.0} leads us to the following Definition.

\begin{definition}\delabel{3.1} Let $(\Cc,x)$ be an $A$-coring with a fixed
grouplike, and let $B=A^{{\rm co}\Cc}$. We
call $(\Cc,x)$ a Galois coring if the canonical coring morphism
${\rm can}:\ \Dd=A\ot_B A\to \Cc$, ${\rm can}(a\ot_B b)=axb$
is an isomorphism.
\end{definition}

Let $i:\ B\to A$ be a ring morphism. If
$x\in G(\Cc)^B$, then we can define a functor
$$\Gamma:\ \Mm^\Dd\to\Mm^\Cc,~~\Gamma(M,\rho)=(M,\widetilde{\rho})$$
with $\widetilde{\rho}(m)=m_{[0]}\ot_A xm_{[1]}\in M\ot_A\Cc$ if
$\rho(m)=m_{[0]}\ot_B m_{[1]}\in M\ot_B A$. It is easy to see that
$\Gamma\circ K=F$, and therefore we have the following result.

\begin{proposition}\prlabel{3.2}
Let $(\Cc,x)$ be a Galois $A$-coring. Then $\Gamma$ is an isomorphism
of categories. Consequently $R$ (resp. $K$) is fully faithful
if and only if $G$ (resp. $F$) is fully faithful.
\end{proposition}

Let us now give some alternative characterizations of Galois corings;
for the proof, we refer to \cite[3.6]{W3}.

\begin{proposition}\prlabel{3.3}
Let $(\Cc,x)$ be an $A$-coring with fixed grouplike element, and
$B=A^{{\rm co}\Cc}$. The following assertions are equivalent.
\begin{enumerate}
\item $(\Cc,x)$ is Galois;
\item if $(M,\rho)\in \Mm^{\Cc}$ is such that $\rho:\ M\to M\ot_A\Cc$
is a coretraction, then the evaluation map
$$\varphi_M:\ \Hom^\Cc(A,M)\ot_BA\to M,~~\varphi_M(f\ot_B m)=f(m)$$
is an isomorphism;
\item $\varphi_\Cc$ is an isomorphism.
\end{enumerate}
\end{proposition}

From \thref{2.4} and \prref{3.2}, we immediately obtain the following
result.

\begin{theorem}\thlabel{3.4}
Let $(\Cc,x)$ be a Galois $A$-coring, and put $B=A^{{\rm
co}\Cc}$. Then the following assertions are equivalent.
\begin{enumerate}
\item $(F,G)$ and $(F',G')$ are equivalences of categories;
\item the functors $F$ and $F'$ are fully faithful;
\item $i:\ B\to A$ is pure in ${}_B\Mm$ and $\Mm_B$.
\end{enumerate}
\end{theorem}

\begin{remark}\relabel{3.5}
Let us make some remarks about terminology.
In the literature, there is an inconsistency in the use of the
term ``Galois". An alternative definition is to require that $(\Cc,x)$
satisfies the equivalent definitions of \thref{3.4}, so that 
$(F,G)$ and $(F',G')$ are category equivalences.
In \seref{5},
we will discuss special cases that have appeared in the literature
before. In some cases, there is an agreement with \deref{3.1}
(see e.g. \cite{Brzezinski99}, \cite{Schneider90}), while in other
cases, category equivalence is required (see e.g. \cite{CS}, \cite{DI}).

In the particular situation where $\Cc=A\ot H$, as in \exref{1.6},
the property that $(F,G)$ is an equivalence (resp. $G$ is fully
faithful)
is called the
Strong (resp. Weak) Structure Theorem (see \cite{DT}).
Let $i:\ B\to A$ be a ring morphism, and $\Dd$ the canonical
coring. In the situation where $A$ and $B$ are commutative,
$i$ is called a descent morphism (resp. an effective descent
morphism) if $K$ is fully faithful (resp. $(K,R)$ is an equivalence).
In the general situation,
$(K,R)$ is an equivalence if and only if the
functor $\bullet\ot_B A:\ \Mm_B\to \Mm_A$ is comonadic
(see e.g. \cite[Ch. 4]{Borceux}).
\end{remark}

Let us next look at the case where $A$ is flat as a left $B$-module.
Wisbauer \cite{W3} calls the following two results the
{\sl Galois coring Structure Theorem}.

\begin{proposition}\prlabel{3.7}
Let $(\Cc,x)$ be an $A$-coring with fixed grouplike element, and
$B=A^{{\rm co}\Cc}$. Then the following statements are equivalent.
\begin{enumerate}
\item $(\Cc,x)$ is Galois and $A$ is flat as a left $B$-module;
\item $G$ is fully faithful
and $A$ is flat as a left $B$-module;
\item $\Cc$ is flat as a left $A$-module, and $A$ is a generator
in $\Mm^\Cc$.
\end{enumerate}
\end{proposition}

\begin{proof}
$\ul{1)\Rightarrow 2)}$ follows from Propositions \ref{pr:2.2}
and \ref{pr:3.2}. $\ul{2)\Rightarrow 1)}$ follows from
\prref{3.0}. For the proof of $\ul{1)\Leftrightarrow 3)}$,
we refer to \cite[3.8]{W3}.
\end{proof}

\begin{proposition}\prlabel{3.8}
Let $(\Cc,x)$ be an $A$-coring with fixed grouplike element, and
$B=A^{{\rm co}\Cc}$. Then the following statements are equivalent.
\begin{enumerate}
\item $(\Cc,x)$ is Galois and $A$ is faithfully flat as a left $B$-module;
\item $(F,G)$ is an equivalence
and $A$ is flat as a left $B$-module;
\item $\Cc$ is flat as a left $A$-module, and $A$ is a projective generator
in $\Mm^\Cc$.
\end{enumerate}
\end{proposition}

\begin{proof}
The equivalence of 1) and 2) follows from Propositions \ref{pr:2.5}
and \ref{pr:3.2}. For the remaining equivalence, we refer to
\cite{W3}.
\end{proof} 

A right $\Cc$-comodule $N$ is called semisimple (resp. simple)
in $\Mm^\Cc$
if every $\Cc$-monomorphism $U\to N$ is a coretraction
(resp. an isomorphism).
Similar definitions apply to left $\Cc$-comodules
and $(\Cc,\Cc)$-bicomodules. $\Cc$ is said to be right (left)
semisimple if it is semisimple as a right (left) $\Cc$-comodule.
$\Cc$ is called a simple coring if it is simple as a
$(\Cc,\Cc)$-bicomodule.
For the proof of the following result, we refer to \cite{KGL}.

\begin{proposition}\prlabel{3.9}
For an $A$-coring $\Cc$, the following assertions are equivalent:
\begin{enumerate}
\item $\Cc$ is right semisimple;
\item $\Cc$ is projective as a left $A$-module and $\Cc$ is
semisimple as a left $\*C$-module;
\item $\Cc$ is projective as a right $A$-module and $\Cc$ is
semisimple as a right $\Cc^*$-module;
\item $\Cc$ is left semisimple.
\end{enumerate}
\end{proposition}

The connection to Galois corings is the following, due to Wisbauer
\cite[3.12]{W3}:

\begin{proposition}\prlabel{3.10}
For an $A$-coring with a fixed grouplike element $(\Cc,x)$, the following
assertions are equivalent:
\begin{enumerate}
\item $\Cc$ is a simple and left (or right) semisimple coring;
\item $(\Cc,x)$ is Galois and $\End^\Cc(A)$ is simple and left
semisimple;
\item $(\Cc,x)$ is Galois and $B$ is a simple left semisimple
subring of $A$;
\item $\Cc$ is flat as a right $A$-module, $A$ is a projective
generator in ${}^\Cc\Mm$, and $\End^\Cc(A)$ is simple and left
semisimple (the left $\Cc$-coaction on $A$ being given by
$\rho^l(a)=ax$).
\end{enumerate}
\end{proposition}

\section{Galois corings and Morita theory}\selabel{4}
Let $(\Cc,x)$ be a coring with a fixed grouplike element,
$B=A^{{\rm co}\Cc}$, and $\Dd=A\ot_BA$. We can consider the
left dual of the map can:
$${}^*\can:\ {}^*\Cc\to {}^*\Dd\cong {}_B\End(A)^{\rm op},~~
{}^*\can(f)(a)=f(xa).$$
The following result is obvious.

\begin{proposition}\prlabel{4.1}
If $(\Cc,x)$ is Galois, then ${}^*\can$ is an isomorphism.
The converse property holds if $\Cc$ and $A$ are finitely generated
projective, respectively as a left $A$-module, and a left $B$-module.
\end{proposition}

Let $Q=\{q\in \*C~|~c_{(1)}q(c_{(2)})=q(c)x,~{\rm for~all~}c\in \Cc\}$.
A straightforward computation shows that $Q$ is a $(\*C,B)$-bimodule.
Also $A$ is a left $(B,\*C)$-bimodule; the right $\*C$-action is
induced by the right $\Cc$-coaction: $a\cdot f=f(xa)$. Now
consider the maps
$$\tau:\ A\ot_{\*C} Q\to B,~~\tau(a\ot_{\*C} q)=q(xa);$$
$$\mu:\ Q\ot_B A\to \*C,~~\mu(q\ot_B a)=q\#i(a).$$
With this notation, we have the following property (see \cite{CVW}).

\begin{proposition}\prlabel{4.2}
$(B,\*C,A,Q,\tau,\mu)$ is a Morita context.
\end{proposition}

Properties of this Morita context are studied in \cite{Abu1},
\cite{Abu2}, \cite{CVW} and \cite{CVW2}. It generalizes (and unifies)
Morita contexts discussed in \cite{BDR}, \cite{CS}, \cite{CF},
\cite{CFM} and \cite{D}. We recall the following properties from
\cite{CVW} and \cite{CVW2}.

\begin{proposition}\prlabel{4.3} \cite[Th. 3.3 and Cor. 3.4]{CVW}
If $\tau$ is surjective, then $M^{{\rm co}\Cc}=
M^{\*C}=\{m\in M~|~m\cdot f=mf(x),~{\rm for~all~}f\in\*C\}$,
for all $M\in \Mm^{\Cc}$.\\
The following assertions are equivalent:
\begin{enumerate}
\item $\tau$ is surjective;
\item there exists $q\in Q$ such that $q(x)=1$;
\item for every $M\in \Mm^{\*C}$, the map
$$\omega_M:\ M\ot_{\*C}Q\to M^{\*C},~~\omega_M(m\ot_{\*C}q)=m\cdot q$$
is bijective.
\end{enumerate}
\end{proposition}

\begin{proposition}\prlabel{4.4} \cite{CVW2}, \cite[Th. 3.5]{CVW}
The following assertions are equivalent:
\begin{enumerate}
\item $\mu$ is surjective;
\item $\Cc$ is finitely generated and projective
as a left $A$-module and $G$ is fully faithful.
\end{enumerate}
\end{proposition}

As an application of \prref{4.3}, we have the following result.

\begin{corollary}\colabel{4.5}
Assume that $\Cc$ is finitely generated projective as a left
$A$-module. Consider the adjoint pair $(F=\bullet\ot_B A,
G=(\bullet)^{{\rm co}\Cc})$, and the functors
$\widetilde{F}=\bullet\ot_B A$ and $\widetilde{G}=\bullet\ot_{\*C} Q$
coming from the Morita context of \prref{4.2}. Then 
$F\cong \widetilde{F}$ and $G\cong \widetilde{G}$ if $\tau$ is surjective.
\end{corollary}

\begin{proof}
Take $N\in \Mm_B$. $F(N)$ corresponds to $\widetilde{F}$ under the
isomorphism $\Mm^{\Cc}\cong \Mm_{\*C}$. If $\tau$ is surjective, then
it follows from \prref{4.3} that $\omega:\ \widetilde{G}\to G$
is bijective.
\end{proof}

Let us now compute the Morita context associated to the canonical
coring.

\begin{proposition}\prlabel{4.6}
Let $i:\ B\to A$ be a ring morphism, and assume that $i$ is pure
as a morphism of left $B$-modules. Then the Morita context associated
to the canonical coring $(\Dd=A\ot_BA,1\ot_B 1)$ is the
Morita context $(B, {}_B\End(A)^{\rm op}, A,$ $ {}_B\Hom(A,B),\varphi,\psi)$
associated to $A$ as a left $B$-module
(see \cite[II.4]{Ba}).
\end{proposition}

\begin{proof}
From \prref{2.1}, it follows that
$$A^{{\rm co}\Dd}=\{b\in A~|~b\ot_B1=1\ot_B b\}=B.$$
Take $q\in Q\subset {}_A\Hom(A\ot_BA,A)$ and the
corresponding $\widetilde{q}\in {}_B\Hom(A,A)$, given by
$\tilde{q}(a)=q(1\ot_B a)$. Then
$$q(a'\ot_B a)(1\ot_B 1)=(a'\ot_B1)q(1\ot_B a).$$
Taking $a'=1$, we find
$$\tilde{q}(a)\ot_B 1=1\ot_B \tilde{q}(a)$$
hence $\tilde{q}(a)\in B$, and it follows that
$Q\subset {}_B\Hom(A,B)$. The converse inclusion is proved in
a similar way. A straightforward verification shows that
$\varphi=\tau$ and $\psi=\mu$.
\end{proof}

Recall that the context associated to the left $A$-module $B$
is strict if and only if $A$ is a left $B$-progenerator.
We are now ready to prove the following result. In \seref{5},
we will see that it is a generalization of \cite[Th. 9.3 and 9.6]{CS}.

\begin{theorem}\thlabel{4.7}
Let $(\Cc,x)$ be a coring with fixed grouplike element, and assume
that $\Cc$ is a left $A$-progenerator. We take a subring
$B'$ of $B=A^{{\rm co}\Cc}$, and consider the map
$$\can':\ \Dd'=A\ot_{B'}A\to \Cc,~~\can'(a\ot_{B'}a')=axa'$$
Then the following statements are equivalent:
\begin{enumerate}
\item \begin{itemize}
\item $\can'$ is an isomorphism;
\item $A$ is faithfully flat as a left $B'$-module.
\end{itemize}
\item \begin{itemize}
\item ${}^*\can'$ is an isomorphism;
\item $A$ is a left $B'$-progenerator.
\end{itemize}
\item \begin{itemize}
\item $B=B'$;
\item the Morita context $(B,\*C,A,Q,\tau,\mu)$ is strict.
\end{itemize}
\item \begin{itemize}
\item $B=B'$;
\item $(F,G)$ is an equivalence of categories.
\end{itemize}
\end{enumerate}
\end{theorem}

\begin{proof}
$\ul{1)\Leftrightarrow 2)}.$ Obviously ${}^*\can'$ is an isomorphism
if $\can'$ is an isomorphism, and the converse holds if
$\Cc$ is a left $A$-progenerator and $A$ is a left $B'$-progenerator.
If $\can'$ is an isomorphism, then $A\ot_{B'}A=\Dd'\cong
\Cc$ is a left $A$-progenerator, hence $A$ is a left $B'$-progenerator.\\
$\ul{1)\Rightarrow 3)}.$ 
Since $A$ is faithfully flat as a left
$B'$-module, $A^{{\rm co}\Dd'}=B'$. Since $\can'$ is an isomorphism,
it follows that $B= A^{{\rm co}\Cc}=A^{{\rm co}\Dd'}=B'$.
Then $\can=\can'$ is an isomorphism, hence the Morita contexts
associated to $(\Cc,x)$ and $(\Dd,1\ot_B 1)$ are isomorphic.
From the equivalence of 1) and 2), we know that $A$ is a left $B$-progenerator,
so the context associated to $(\Dd,1\ot_B 1)$ is strict, 
see the remark preceding \thref{4.7}.
Therefore the Morita context
$(B,\*C,A,Q,\tau,\mu)$ associated to $(\Cc,x)$ is also strict.\\
$\ul{3)\Rightarrow 1)}.$ and $\ul{3)\Rightarrow 4)}.$
If $(B,\*C,A,Q,\tau,\mu)$ is strict, then
$A$ is a left $B$-progenerator, and a fortiori faithfully flat as
a left $B$-module. $\tau$ is surjective, so it follows from 
\coref{4.5} that $F\cong \widetilde{F}$ and 
$G\cong \widetilde{G}$. $(\widetilde{F},\widetilde{G})$ is an
equivalence, so $(F,G)$ is also an equivalence. Then
$(\Cc,x)$ is
Galois by \prref{3.0}.\\
$\ul{4)\Rightarrow 1)}.$ $\can$ is an isomorphism, by
 \prref{3.0}, and we have already seen that this implies that
$A$ is a left $B$-progenerator, so $A$ is faithfully flat as a left
$B$-module.
\end{proof}

\section{Application to particular cases}\selabel{5}
\subsection{Coalgebra Galois extensions}\selabel{5.1}
From \cite{BrzezinskiH}, we recall the following Definition.

\begin{definition}\delabel{5.1}
Let $i:\
B\to A$ be a morphism of $k$-algebras, and $C$ a $k$-coalgebra.
$A$ is called a $C$-Galois extension
of $B$ if the following conditions hold:
\begin{enumerate}
\item $A$ is a right $C$-comodule;
\item ${\rm can}:\ A\ot_B A\to A\ot C$, ${\rm can}(a\ot_B a')=aa'_{[0]}
\ot_B a'_{[1]}$ is an isomorphism;
\item $B=\{a\in A~|~\rho(a)=a\rho(1)\}$.
\end{enumerate}
\end{definition}

\begin{proposition}\prlabel{5.2}
Let $i:\
B\to A$ be a morphism of $k$-algebras, and $C$ a $k$-coalgebra.
$A$ is called a $C$-Galois extension
of $B$ if and only if there exists a right-right entwining
structure $(A,C,\psi)$ and $x\in G(A\ot C)$
such that $A^{{\rm co}A\ot C}=B$ and $(A\ot C,x)$ is a Galois
coring.
\end{proposition}

\begin{proof}
Let $(A,C,\psi)$ be an entwining structure. We have seen in
\exref{1.7} that $\Cc=A\ot C$ is a coring. Given a grouplike
element $x=\sum_i a_i\ot c_i$, we have a right $\Cc$-coaction
on $A$, hence $A$ is an entwined module (see \exref{1.7}),
and therefore a $C$-comodule. The $C$-coaction is given by
the formula
$$\rho(a)=a_{[0]}\ot a_{[1]}=\sum_i a_ia_\psi\ot c_i^\psi.$$
Then the conditions of \deref{5.1} are satisfied, and $A$ is
$C$-Galois extension of $B$.\\
Conversely, let $A$ be a $C$-Galois extension of $B$.
${\rm can}$ is bijective, so the coring structure on $A\ot_B A$
induces a coring structure on $A\ot C$. We will show that this coring
structure comes from an entwining structure $(A,C,\psi)$.\\
It is clear that the natural left $A$-module structure on $A\ot C$
makes ${\rm can}$ into a left $A$-linear map.
The right $A$-module structure on $A\ot C$ induced by can is given by
$$(b\ot c)a={\rm can}({\rm can}^{-1}(b\ot c)a).$$
Since ${\rm can}^{-1}(1_{[0]}\ot 1_{[1]})=1\ot_B 1$, we have
\begin{equation}\eqlabel{5.2.1}
(1_{[0]}\ot 1_{[1]})a={\rm can}(1\ot a)=a_{[0]}\ot a_{[1]}.
\end{equation}
The comultiplication $\Delta$ on $A\ot C$ is given by
$$\Delta(a\ot c)=({\rm can}\ot_A{\rm can})\Delta_{\Cc}({\rm can}^{-1}
(a\ot c))\in (A\ot C)\ot_A (A\ot C),$$
for all $a\in A$ and $c\in C$. can is bijective, so we can find
$a_i,b_i\in A$ such that
$${\rm can}(\sum_i a_i\ot_B b_i)=\sum_i a_ib_{i[0]}\ot_B b_{i[1]}=a\ot c,$$
and we compute that
\begin{eqnarray*}
\Delta(a\ot c)&=&
({\rm can}\ot_A{\rm can})\Delta_{\Cc}(\sum_i a_i\ot_B b_i)\\
&=&\sum_i {\rm can}(a_i\ot_B 1)\ot_A {\rm can}(1\ot_B b_i)\\
&=& \sum_i (a_i1_{[0]}\ot 1_{[1]})\ot_A (b_{i[0]}\ot b_{i[1]})\\
&=& \sum_i (a_i1_{[0]}\ot 1_{[1]})b_{i[0]}\ot_A (1\ot b_{i[1]})\\
{\rm \equref{5.2.1}}~~~~
&=& \sum_i (a_ib_{i[0]}\ot b_{i[1]})\ot_A (1\ot b_{i[2]})\\
&=& (a\ot c_{(1)})\ot_A (1\ot c_{(2)}).
\end{eqnarray*}
Finally
\begin{eqnarray*}
\varepsilon_Cc(a\ot c)&=& \varepsilon_{\Cc}(\sum_i a_i\ot b_i)=
\sum_i a_i b_i\\
&=& \sum_i a_ib_{i[0]}\varepsilon_C(b_{i[1]})
=a\varepsilon_C(c).
\end{eqnarray*}
Now define $\psi:\ C\ot A\to A\ot C$ by $\psi(c\ot a)=(1_A\ot c)a$.
It follows from \cite{Br} that $(A,C,\psi)$ is an entwining structure.
\end{proof}

Let $(A,C,\psi)$ be an entwining structure, and consider $g\in C$
grouplike. Then $x=1_A\ot g$ is a grouplike element of $A\ot C$.
Let us first describe the Morita context from the previous Section.\\
Observe that
${\*C}={}_A\Hom(A\ot C,A)\cong \Hom(C,A)$
as a $k$-module. The ring structure on ${\*C}$ induces a
$k$-algebra structure on $\Hom(C,A)$, and this $k$-algebra is
denoted $\#(C,A)$. The product is given by the formula
\begin{equation}\eqlabel{3.1}
(f\# g)(c)=f(c_{(2)})_{\psi}g(c_{(1)}^{\psi}).
\end{equation}
We have a natural algebra homomorphism $i:\ A\to \#(C,A)$, 
$i(a)(c)=\varepsilon_C(c)a$,
and we have, for all $a\in A$ and $f:\ C\to A$:
\begin{equation}\eqlabel{3.3}
(i(a)\# f)(c)=a_{\psi}f(c^{\psi})~~{\rm and}~~
(f\#i(a))(c)=f(c)a.
\end{equation}
$\Hom(C,A)$ will denote the $k$-algebra with the usual convolution
product, that is
\begin{equation}\eqlabel{3.4}
(f* g)(c)=f(c_{(1)})g(c_{(2)}).
\end{equation}
The ring of coinvariants is
\begin{equation}\eqlabel{3.6}
B=A^{{\rm co}C}=\{b\in A~|~b_{\psi}\ot g^{\psi}=b\ot g\},
\end{equation}
and the bimodule $Q$ is naturally isomorphic
to
$$Q=\{q\in \#(C,A)~|~q(c_{(2)})_{\psi}\ot c_{(1)}^{\psi}=q(c)\ot g\}.$$
We have maps
$$\mu:\ Q\ot_{B'} A\to\#(C,A),~~\mu(q\ot_{B} a)(c)=q(c)a,$$
$$\tau:\ A\ot_{\#(C,A)} Q\to B,~~\tau(a\ot q)=a_{\psi}q(x^{\psi}),$$
and $(B,\#(C,A),A,Q,\tau,\mu)$ is a Morita context.

\begin{proposition}\prlabel{5.3} \cite[Prop. 4.3]{CVW}
Assume that $\lambda:\ C\to A$ is convolution invertible, with
convolution inverse $\lambda^{-1}$. Then the following assertions are
equivalent:
\begin{enumerate}
\item $\lambda\in Q$;
\item for all $c\in C$, we have
\begin{equation}\eqlabel{5.3.1}
\lambda^{-1}(c_{(1)})\lambda(c_{(3)})_{\psi}\ot c_{(2)}^{\psi}=
\varepsilon(c)1_A\ot g;
\end{equation}
\item for all $c\in C$, we have
\begin{equation}\eqlabel{5.3.2}
\lambda^{-1}(c_{(1)})\ot c_{(2)}=\lambda^{-1}(c)_{\psi}\ot g^{\psi}.
\end{equation}
\end{enumerate}
Notice that condition 3) means that $\lambda^{-1}$ is right
$C$-colinear. If such a $\lambda\in Q$ exists, then we call
$(A,C,\psi,g)$ cleft.
\end{proposition}

\begin{proposition}\prlabel{5.4} \cite[Prop. 4.4]{CVW}
Assume that $(A,C,\psi,g)$ is a cleft entwining structure. Then the
map $\tau$ in the associated Morita context is surjective.
\end{proposition}

We say that the entwining structure $(A,C,\psi,g)$ satisfies the {\sl right
normal basis property} if there exists a left $B$-linear and
right $C$-colinear isomorphism
$B\ot C\to A$. The following is one of the main results in \cite{CVW}.
As before, we consider the functor $F=\bullet\ot_B A:\ \Mm_B\to
\Mm(\psi)_A^C$ and its right adjoint $G=(\bullet)^{{\rm co}C}$.

\begin{theorem}\thlabel{5.5} \cite[Theorem 4.5]{CVW}
Let $(A,C,\psi,g)$ be an entwining structure with a fixed grouplike element.
The following assertions are equivalent:
\begin{enumerate}
\item $(A,C,\psi,g)$ is cleft;
\item $(F,G)$ is a category equivalence and
$(A,C,\psi,g)$ satisfies the right normal basis property;
\item $(A,C,\psi,g)$ is Galois, and satisfies the right normal basis property;
\item the map ${}^*{\rm can}:\ \#(C,A)\to \End_B(A)^{\rm op}$ is
bijective and $(A,C,\psi,g)$ satisfies the
right normal basis property.
\end{enumerate}
\end{theorem}

\subsection{Hopf-Galois extensions}\selabel{5.2}
Let $H$ be a Hopf algebra over a commutative ring $k$ with
bijective antipode, and $A$ a right $H$-comodule algebra
(cf. \exref{1.6}). Then $\Cc=A\ot H$ is an $A$-coring,
and $1_A\ot 1_H\in G(\Cc)$. Let $B=A^{{\rm co}H}$.
The canonical map is now the following:
$$\can:\ A\ot_B A\to A\ot H,~~\can(a'\ot_B a)=a'a_{[0]}\ot a_{[1]}$$

\begin{definition}\delabel{5.6} (see e.g. \cite[Def. 1.1]{DT})
$A$ is a Hopf-Galois extension of $B$ if and only if
$\can$ is an isomorphism.
\end{definition}

Obviously $A$ is a Hopf-Galois extension of $B$ if and only if
$(A\ot H,1_A\ot 1_H)$ is a Galois coring.\\

Assume now that $H$ is a progenerator as a $k$-module, i.e.
$H$ is finitely generated, faithful, and projective as a $k$-module.
Then $A\ot H$ is a left $A$-progenerator, so we can apply the
results of \seref{4}. We will show that we recover results from
\cite{CS}.
To this end, we will describe the Morita context associated to
$(A\ot H,1_A\ot 1_H)$.\\
First we compute $\*C$. We have already seen in \seref{5.1}
that $\*C\cong \#(H,A)$. As a module, $\Hom(H,A)\cong H^*\ot A$,
since $H$ is finitely generated and projective. The multiplication on
$\#(H,A)$ can be transported into a multiplication on $H^*\ot A$,
giving us a $k$-algebra denoted by $H^*\# A$.
A straightforward computation shows that this multiplication
is given by the following formula. $H^*$ is a coalgebra, since
$H$ is finitely generated projective, and $H^*$ acts on $A$
from the left: $h^*\leftact a=\lan h^*,a_{[1]}\ran a_{[0]}$.
Then we can compute that
\begin{equation}\eqlabel{5.7.1}
(h^*\# a)(k^*\# b)=(k^*_{(1)}* h^*)\# (k^*_{(2)}\leftact a_{[0]})b.
\end{equation}
Consider the map $\can':\ A\ot A\to A\ot H$; its dual
${}^*\can':\ H^*\# A\to \End(A)^{\rm op}$ is given by
\begin{equation}\eqlabel{5.7.2}
{}^*\can'(h^*\# a)(b)=(h^*\leftact b)a.
\end{equation}
Take $y=\sum_i h_i^*\# a_i\in H^*\# A$. $y\in Q$ if and only if
$$\sum_i \lan h_i^*, h_{(2)}\ran a_{i[0]}\ot h_{(1)}a_{i[1]}=
\sum_i\lan h_i^*,h\ran a_i\ot 1,$$
for $h\in H$. Since $H$ is finitely generated and projective, this is also 
equivalent to
$$\sum_i \lan h_i^*, h_{(2)}\ran a_{i[0]}\lan h^*, h_{(1)}a_{i[1]}\ran=
\sum_i\lan h_i^*,h\ran a_i\lan h^*, 1\ran,$$
for all $h\in H$ and $h^*\in H^*$. The left hand side equals
$$\sum_i \lan h^*_{(1)}*h_i^*,h\ran \lan h^*_{(2)},a_{i[1]}\ran a_{i[0]},$$
so we find that $y\in Q$ if and only if
$$\sum_i (h^*_{(1)}*h_i^*)\# (h_{(2)}^*\leftact a_i)=
\lan h^*, 1\ran \sum_i h_i^*\#a_i,$$
or
$$y(h^*\#1)= \lan h^*, 1\ran y.$$
for all $h^*\in H^*$. Thus
\begin{equation}\eqlabel{5.7.3}
Q=\{y\in H^*\# A~|~y(h^*\#1)= \lan h^*, 1\ran y,~{\rm for~all~}h^*\in H^*\}.
\end{equation}
Elementary computations show that the maps $\mu$ and $\tau$ from the
Morita context are the following:
$$\tau:\ A\ot_{H^*\# A}Q\to B,~~
\tau(a\ot y)={}^*\can(y)(a);$$
$$\mu:\ Q\ot_B A\to H^*\# A,~~\mu(y\ot a)=y(\varepsilon_C\# a),$$
where $B=A^{{\rm co}H}$, as usual.
\thref{4.7} now takes the following form.

\begin{proposition}\prlabel{5.7}
Let $H$ be a $k$-progenerator Hopf algebra over a commutative ring $k$,
and $A$ a right $H$-comodule algebra. Then the following statements
are equivalent (with notation as above):
\begin{enumerate}
\item
\begin{itemize}
\item $\can':\ A\ot A\to A\ot H,~~\can'(a'\ot a)=a'a_{[0]}\ot a_{[1]}$
is bijective;
\item $A$ is faithfully flat as a $k$-module.
\end{itemize}
\item
\begin{itemize}
\item ${}^*\can':\ H^*\# A\to \End(A)^{\rm op},~~{}^*\can(h^*\# a)(b)=
(h^*\leftact b)a$ is an isomorphism;
\item $A$ is a $k$-progenerator.
\end{itemize}
\item
\begin{itemize}
\item $A^{{\rm co}H}=k$;
\item the Morita context $(k,H^*\# A,A,Q,\tau,\mu)$ is strict.
\end{itemize}
\item
\begin{itemize}
\item $A^{{\rm co}H}=k$;
\item the adjoint pair of functors 
$(F=\bullet\ot A,G=(\bullet)^{{\rm co}H})$
is an equivalence between the categories $\Mm_k$ and $\Mm_A^H$.
\end{itemize}
\end{enumerate}
\end{proposition}

\cite[Theorems 9.3 and 9.6]{CS} follow from \prref{5.7}.

\subsection{Classical Galois Theory}\selabel{5.3}
As in \exref{1.4}, let $G$ be a finite group, and $A$ a $G$-module
algebra. We have seen that $\Cc= A\ot (kG)^*=\oplus_{\sigma\in G}
Av_\sigma$ is an $A$-coring. $\sum_{\sigma}v_\sigma$ is
a grouplike element. Since $(kG)^*$ is finitely generated and
projective, we can apply \prref{5.7}.
We have
$$\can':\ A\ot A\to \oplus_{\sigma\in G} Av_\sigma,~~
\can'(a\ot b)=\sum_\sigma a\sigma(b)v_\sigma.$$
$$\*C=\oplus_\sigma u_\sigma A,$$
with multiplication
$$(u_\sigma a)(u_\tau b)=u_{\tau\sigma}\tau(a)b,$$
and
$${}^*\can':\ \oplus_\sigma u_\sigma A\to \End(A)^{\rm op},~~
{}^*\can'(u_\sigma a)(b)=\sigma(b)a.$$
We also have
$$Q=\{\sum_\sigma u_\sigma \sigma(a)~|~a\in A\}\cong A,$$
which is not surprising since $(kG)^*$ is a Frobenius Hopf
algebra (see \cite{CVW} and \cite{CFM}).
If $A^G=k$, then we have a Morita context
$(k,\*C,A,A,\tau,\mu)$, where the connecting maps are the
following:
$$\tau:\ A\ot_{\*C} A\to k,~~\tau(a\ot b)=\sum_{\sigma}\sigma(ab);$$
$$\mu:\ A\ot A\to \*C,~~\mu(a\ot b)=\sum_\sigma u_\sigma\sigma(b)a.$$
\prref{5.7} now takes the following form (compare to \cite[Prop.
III.1.2]{DI}).

\begin{proposition}\prlabel{5.8}
Let $G$ be a finite group, $k$ a commutative ring and $A$
a $G$-module algebra. Then the following statements are equivalent:
\begin{enumerate}
\item
\begin{itemize}
\item $\can'$ is an isomorphism;
\item $A$ is faithfully flat as a $k$-module.
\end{itemize} 
\item
\begin{itemize}
\item ${}^*\can'$ is an isomorphism;
\item $A$ is a $k$-progenerator.
\end{itemize}
\item
\begin{itemize}
\item $A^G=k$;
\item the Morita context $(k,\*C,A,A,\tau,\mu)$ is strict.
\end{itemize}
\item
\begin{itemize}
\item $A^G=k$;
\item the adjoint pair of functors $(F=\bullet\ot A,G=(\bullet)^G)$
is an equivalence between the categories of $k$-modules and
right $A$-modules on which $G$ acts as a group of right $A$-semilinear
automorphisms.
\end{itemize}
\end{enumerate}
\end{proposition}

In the case where $A$ is a commutative $G$-module algebra, we
have some more equivalent conditions.

\begin{proposition}\prlabel{5.8b}
Let $G$ be a finite group, $k$ a commutative ring and $A$
a commutative $G$-module algebra. Then the statements of
\prref{5.8} are equivalent to
\begin{enumerate}
\item[5.]
\begin{itemize}
\item $A^G=k$;
\item for each non-zero idempotent $e\in S$ and $\sigma\neq \tau\in G$,
there exists $a\in A$ such that $\sigma(a)e\neq \tau(a)e$;
\item $A$ is a separable $k$-algebra (i.e. $A$ is projective as
an $A$-bimodule).
\end{itemize}
\item[6.] 
\begin{itemize}
\item $A^G=k$;
\item there exist $x_1,\cdots,x_n,y_1,\cdots y_n\in A$
with
$$\sum_{j=1}^n x_j\sigma(y_j)=\delta_{\sigma,e}$$
for all $\sigma\in G$.
\end{itemize}
\item[7.] 
\begin{itemize}
\item $A^G=k$;
\item for each maximal ideal $m$ of $A$, and for each
$\sigma\neq e\in G$, there exists $x\in S$ such that
$\sigma(x)-x\not\in m$.
\end{itemize}
\end{enumerate}
\end{proposition}

\begin{proof}
We refer to \cite[Th 1.3]{CHR} and \cite[Prop. III.1.2]{DI}.
\end{proof}

If $A=l$ is a field, then the second part of condition 7.
is satisfied. Let $l$ be a finite field extension of a field
$k$, and $G$ the group of $k$-automorphisms of $l$. Then
$l^G=k$ if and only if $l$ is a normal and separable
(in the classical sense) extension of $k$ (see e.g.
\cite[Th. 10.8 and 10.10]{Stewart}). Thus we recover the classical
definition of a Galois field extension.

\subsection{Strongly graded rings}\selabel{5.4}
As in \exref{1.5}, let $G$ be a group, and $A$ a $G$-graded
ring, and $\Cc=\oplus_{\sigma\in G}Au_{\sigma}$. Fix
$\lambda\in G$, and take the grouplike element
$u_\lambda\in G(\Cc)$. Then $M^{{\rm co}\Cc}=M_\lambda$,
for any right $G$-graded $A$-module,
and $B=A^{{\rm co}\Cc}=A_e$. Since $B$ is a direct factor
of $A$, $A$ is flat as a left and right $B$-module,
and $i:\ B\to A$ is pure in $\Mm_B$ and ${}_B\Mm$. Also
$$\can:\ A\ot_B A\to \oplus_{\sigma\in G}Au_{\sigma},~~
\can(a'\ot a)=\sum_{\sigma\in G}a'a_\sigma u_{\lambda\sigma}.$$

\begin{proposition}\prlabel{5.9}
With notation as above, the following assertions are equivalent.
\begin{enumerate}
\item $A$ is strongly $G$-graded, that is, $A_\sigma A_\tau=
A_{\sigma \tau}$, for all $\sigma,\tau\in G$;
\item the pair of adjoint functors 
$(F=\bullet\ot_BA,G=(\bullet)_\lambda)$ is an equivalence 
between $\Mm_B$ and $\Mm_A^G$, the category of $G$-graded
right $A$-modules;
\item $(\Cc,u_{\lambda})$ is a Galois coring.
\end{enumerate}
In this case $A$ is faithfully flat as a left (or right) $B$-module.
\end{proposition}

\begin{proof}
$\ul{1)\Rightarrow 2)}$ is a well-known fact from graded ring
theory. We sketch a proof for completeness sake. The unit of
the adjunction between $\Mm_B$ and $\Mm_A^G$ is given by
$$\eta_N:\ N\to (N\ot_B A)_\lambda,~~\eta_N(n)=n\ot_B 1_A.$$
$\eta_N$ is always bijective, even if $A$ is not strongly graded.
Let us show that the counit maps $\zeta_M:\ M_\lambda\ot_B A\to
M$, $\zeta_M(m\ot_B a)=ma$ are surjective. For each $\sigma\in G$,
we can find $a_i\in A_{\sigma^{-1}}$ and $a'_i\in A_{\sigma}$ such
that $\sum_i a_ia'_i=1$. Take $m\in M_\tau$ and
put $\sigma=\lambda_{-1}\tau$. Then
$m=\zeta_M(\sum_i ma_i\ot_B a'_i)$, and $\zeta_M$ is surjective.\\
If $m_j\in M_\lambda$ and $c_j\in A$ are such that
$\sum_j m_jc_j=0$, then for each $\sigma\in G$, we have
$$\sum_j m_j\ot c_{j\sigma}=
\sum_{i,j} m_j\ot c_{j\sigma}a_ia'_i=
\sum_{i,j} m_jc_{j\sigma}a_i\ot a'_i=0.$$
hence $\sum_j m_j\ot c_{j}=\sum_{\sigma\in G}=
\sum_j m_j\ot c_{j\sigma}=0$, so $\zeta_M$ is also injective.\\
$\ul{2)\Rightarrow 3)}$ follows from \prref{3.0}.\\
$\ul{3)\Rightarrow 1)}$ follows from \thref{3.4} and the
fact that $i:\ B\to A$ is pure in $\Mm_B$ and ${}_B\Mm$.\\
The final statement follows from \prref{2.5} and the
fact that $A$ is flat as a $B$-module. 
\end{proof}

Notice that, in this situation, the fact that $(\Cc, u_\lambda)$
is Galois is independent of the choice of $\lambda$.

\section{A more general approach: comatrix corings}\selabel{6}
Let $\Cc$ be an $A$-coring. We needed a grouplike $x\in \Cc$
in order to make $A$ into a right $\Cc$-comodule. In
\cite{Kaoutit}, the following idea is investigated. A couple
$(\Cc,\Sigma)$, consisting of a coring $\Cc$
and a right $\Cc$-comodule $\Sigma$ which is finitely generated
and projective as a right $A$-module, will be called a coring
with a fixed finite comodule.
Let $T=\End^\Cc(\Sigma)$.
Then we have a pair of adjoint functors
$$F=\bullet\ot_T \Sigma:\ \Mm_T\to \Mm^\Cc~~;~~
G=\Hom^\Cc(\Sigma,\bullet):\ \Mm^\Cc\to \Mm_T,$$
with unit $\nu$ and counit $\zeta$ given by
$$\nu_N:\ N\to \Hom^\Cc(\Sigma,N\ot_T\Sigma),~~
\nu_N(n)(u)=n\ot_Tu;$$
$$\zeta_M:\ \Hom^\Cc(\Sigma,M)\ot_T\Sigma\to M,~~
\zeta_M(f\ot_T u)=f(u).$$
In the situation where $\Sigma=A$, we recover the adjoint pair
discussed at the end of \seref{1}. 
A particular example is the {\sl comatrix coring}, generalizing
the {\sl canonical coring}. Let $A$ and $B$ be rings,
and $\Sigma\in {}_B\Mm_A$, with $\Sigma$ finitely generated
and projective as a right $A$-module. Let
$$\{(e_i,e_i^*)~|~i=1,\cdots,n\}\subset \Sigma\times\Sigma^*$$
be a finite dual basis of $\Sigma$ as a right $A$-module.
$\Dd=\Sigma^*\ot_B\Sigma$ is an $(A,A)$-bimodule, and an
$A$-coring, via
$$\Delta_\Dd(\varphi\ot_B u)=
\sum_i\varphi\ot_B e_i\ot_Ae_i^*\ot_B u~~{\rm and}~~
\varepsilon_\Dd(\varphi\ot_B u)=\varphi(u).$$
Furthermore $\Sigma\in \Mm^\Dd$ and $\Sigma^*\in {}^\Dd\Mm$.
The coactions are given by
$$\rho^r(u)= \sum_i e_i\ot_Ae_i^*\ot_B u~~;~~
\rho^l(\varphi)=\sum_i\varphi\ot_B e_i\ot_Ae_i^*.$$
We also have that ${}^*\Dd\cong {}_B\End(\Sigma)^{\rm op}$.
El Kaoutit and G\'omez Torrecillas proved the following
generalization of the Faithfully Flat Descent Theorem.

\begin{theorem}\thlabel{6.1}
Let $\Sigma\in {}_B\Mm_A$ be finitely generated and projective
as a right $A$-module, and $\Dd=\Sigma^*\ot \Sigma$. Then
the following are equivalent
\begin{enumerate}
\item $\Sigma$ is faithfully flat as a left $B$-module;
\item $\Sigma$ is flat as a left $B$-module and
$(\bullet\ot_B\Sigma,\Hom^\Dd(\Sigma,\bullet))$ is a category
equivalence between $\Mm_B$ and $\Mm^\Dd$.
\end{enumerate}
\end{theorem}

Let $(\Cc,\Sigma)$ be a coring with a fixed finite comodule,
and $T= \End^\Cc(\Sigma)$. We have
an isomorphism $f:\ \Sigma^*\to \Hom^\Cc(\Sigma,\Cc)$
given by
$$f(\varphi)=(\varphi\ot_AI_\Cc)\circ\rho~{\rm and}~
f^{-1}(\phi)=\varepsilon_\Cc\circ\phi,$$
for all $\varphi\in \Sigma^*$ and $\phi\in \Hom^\Cc(\Sigma,\Cc)$.
Consider the map
$$\can=\zeta_\Cc\circ (f\ot_BI_\Sigma):\
\Dd=\Sigma^*\ot_T\Sigma\to  \Hom^\Cc(\Sigma,\Cc)\ot_T\Sigma
\to \Cc.$$
We compute easily that $\can(\varphi\ot_B u)=\varphi(u_{[0]})u_{[1]}$.
$\can$ is a morphism of corings, and $\can$ is an isomorphism
if and only if $\zeta_\Cc$ is an isomorphism.

\begin{definition}\delabel{6.2} \cite[3.4]{Kaoutit}
Let $(\Cc,\Sigma)$ be a coring with a fixed finite comodule, and $T=\End^\Cc(\Sigma)$.
$(\Cc,\Sigma)$ is termed Galois if $\can:\ \Sigma^*\ot_T\Sigma\to \Cc$
is an isomorphism.
\end{definition}

\begin{theorem}\thlabel{6.3} \cite[3.5]{Kaoutit}
If $(\Cc,\Sigma)$ is Galois, and $A$ is faithfully flat as
a left $T$-module, then $(F,G)$ is an equivalence of categories.
\end{theorem}

For further results, we refer to \cite{Kaoutit}.

\begin{center}
{\sc Acknowledgements}
\end{center}
The author thanks George Janelidze for stimulating discussions
about descent theory, and the referee for pointing out that
the proper assumption in \prref{2.3} is right purity of $i$,
instead of left purity.

\end{document}